\documentclass[12pt,a4paper]{article}
\usepackage{amsmath,amssymb,bm,ascmac,bbm,mathtools}
\usepackage[dvipdfmx, usenames]{color}
\usepackage[dvipdfmx]{graphicx}
\usepackage{xcolor}
\usepackage{here}
\usepackage{tensor}
\usepackage{authblk}
\usepackage[hang,small,bf]{caption}
\usepackage[subrefformat=parens]{subcaption}
\usepackage{dcolumn}

\newcommand{\mcal}{\mathcal}
\newcommand{\mbb}{\mathbb}
\newcommand{\mfrak}{\mathfrak}

\newcommand{\rank}{\text{rank}}

\newcommand{\mc}{\mathcal}

\newcommand{\wf}[1]{\widehat{\mfrak{#1}}}
\newcommand{\vecG}{\text{Vec}}
\newcommand{\vect}{\text{Vect}_{\mbb C}}
\newcommand{\fp}{\text{FPdim}}
\newcommand{\tc}{\text{ToricCode}}
\newcommand{\fib}{\text{Fib}}
\newcommand{\ising}{\text{Ising}}

\newcommand{\mods}{\text{mod }}

\begin{document}
\title{Classification of connected étale algebras in multiplicity-free modular fusion categories at rank six}
\author{Ken KIKUCHI, Kah-Sen KAM and Fu-Hsiang Huang\footnote{r10222098@ntu.edu.tw}}
\affil{Department of Physics, National Taiwan University, Taipei 10617, Taiwan}
\date{}
\maketitle

\begin{abstract}
We classify connected \'etale algebras $A$'s in multiplicity-free modular fusion categories (MFCs) $\mathcal{B}$'s at rank six, namely $\text{rank}(\mathcal{B})=6$. There are eight MFCs in total and the result indicates that only $so(5)_2$ has a nontrivial connected \'etale algebra. We briefly mention anyon condensation as it is used to determine the category of right $A$-modules in $so(5)_2$. Finally, we discuss physical applications, specifically proving spontaneous $\mathcal{B}$-symmetry breaking (SSB) of these MFCs. The discussion also includes predicting ground state degeneracies and SSB in massive renormalization group flows from two non-unitary minimal models.
\end{abstract}

\makeatletter
\renewcommand{\theequation}
{\arabic{section}.\arabic{equation}}
\@addtoreset{equation}{section}
\makeatother

\section{Introduction}
In this paper, $\mathcal{B}$ denotes a modular fusion category (MFC) (see \cite{EGNO15} for definitions) and $b_j$'s denotes its elements (simple objects) respectively. (For just fusion categories, $\mathcal{C}$ and $c_j$'s serve the same purpose.) The cardinality of an MFC $\mathcal{B}$ is called rank($\mathcal{B}$), so the index $j=1,\ldots,\text{rank}(\mathcal{B})$. The purpose of this paper is to classify connected \'etale algebras in multiplicity-free MFCs at ranks six \cite{GK94,LPR20,VS22} summarized in AnyonWiki \cite{anyonwiki}.\footnote{MFCs up to rank five have been classified in \cite{GK94,RSW07,BNRW15}.} For classification work at rank lower than six, see \cite{KK23preMFC,KK23rank5}. The main result of this paper is summarized as follows.
\vspace{300pt}

\textbf{Theorem.} \textit{Connected \'etale algebras in multiplicity-free modular fusion categories at rank six are given by}
\begin{table}[H]
\begin{center}
\makebox[1 \textwidth][c]{       
\resizebox{1.2 \textwidth}{!}{\begin{tabular}{c|c|c|c}
    Rank&$\mcal B$&Results&Completely anisotropic?\\\hline\hline
    6&$\vecG_{\mbb Z/6\mbb Z}^\alpha$&\cite{G23}, Table \ref{rank6Z6results}&Yes\\
    &$\vecG_{\mbb Z/2\mbb Z}^{-1}\boxtimes\ising$& Table \ref{rank6BFR601results}&Yes\\
    &$su(3)_2\simeq\fib\boxtimes\vecG_{\mbb Z/3\mbb Z}^1$&\cite{EP09}, Table \ref{rank6su32results}&Yes\\
    &$\text{TriCritIsing}$&\cite{KL02}, Table \ref{rank6tricritisingresults}&Yes\\
    &$su(2)_5\simeq\vecG_{\mbb Z/2\mbb Z}^{-1}\boxtimes psu(2)_5$&\cite{KO01}, Table \ref{rank6su25results}&Yes\\
    &$so(5)_2$&\cite{G23}, Table \ref{rank6so52results}&No\\
    &$\fib\boxtimes psu(2)_5$& Table \ref{rank6fibpsu25results}&Yes\\
    &$psu(2)_{11}$& Table \ref{rank6fibpsu211results}&Yes
\end{tabular}.}}
\end{center}
\caption{Connected étale algebras in multiplicity-free MFC $\mcal B$ at rank six}\label{results}
\end{table}

\textbf{Remark.} Some MFCs are realized by Wess-Zumino-Witten models or minimal models. In those cases, we collectively denote the MFCs sharing the same fusion ring by the realization, e.g., $su(3)_2$ or $\text{TriCritIsing}$. Some MFCs are realized by subcategory of objects invariant under centers. We denote the MFCs by realization with $p$ in their head, e.g., $psu(2)_{11}$.\newline

\textbf{Remark.} The classification problem has been actively studied since its inception \cite{DMNO10,DMO11}. Especially, many results for MFCs $\mc C(\mfrak g,k)$'s realized by $\wf g_k$ WZW models are known. (In this context, connected étale algebras are called quantum subgroups \cite{O00}.) For instance, connected étale algebras were classified in \cite{KO01} (for $\wf{su}(2)_k$), in \cite{EP09} (for $\wf{su}(3)_k$), in \cite{CEM23} (for $\wf{su}(4)_k$), and many more in \cite{G23}. Also, connected étale algebras in minimal models with $c<1$ were classified in \cite{KL02}. When available, our results are consistent with them; $\mc C(A_5,1),\mc C(A_2,2),\mc C(A_1,5)$ and the MFC of the tricritical Ising model are known to be completely anisotropic, and $\mc C(C_2,2)$ has two connected étale algebras $A\cong1$ and the $\mbb Z/2\mbb Z$ algebra corresponding to the center of $C_2$. These results are specialized to unitary MFCs, but we classify connected étale algebras in all MFCs including non-unitary ones. \newline

While the classification problem was formulated in mathematics, it has many applications in physics including, say, classification of modular invariants, gauging of categorical symmetries, and anyon condensation. Here, before laying out the details in sections below, we would like to give a prelude to the concept of anyon condensation \cite{BSS02,BSS02',BS08} since it actually plays a crucial role in the classification work of this paper and other \cite{KK23rank5}, in  particular the identification of the correct category of right $A$-modules associated with a non-trivial \'etale algebra. In this paper, we think of the simple objects of each modular fusion categories (MFC) as topological charges or anyon types. Therefore, we will say each MFC in our examples is labeled by some finite set of topological charges or charge sectors or anyonic charges, and these terms are equivalent and will be used interchangeably.

First of all, by anyon condensation, we mean condensation of bosonic anyons, according to \cite{BS08}. (Physically natural conditions force condensable anyons be connected étale \cite{K13}.) The anyon condensation induces a transition between two topologically ordered phase, both described by unitary modular tensor categories (UMTC). That is, this condensation constitutes a symmetry breaking scheme, although now the symmetry is described by a tensor categories or quantum groups. However, this transition is different from the counterpart arising from spontaneous symmetry breaking in which a group is broken to a subgroup: anyon condensation always connect gapped phases with distinct topological orders and there is no Goldstone mode \cite{B17}. In order for an anyon \textit{a} to be a boson, it must satisfy two conditions: (1) trivial spin factor $\theta_a=1$, recalling that $\theta_a$ is related to the conformal dimension $h_a$ by $\theta_a=e^{2\pi ih_a}$, which in turn implies $h_a\in \mathbb{Z}$. (2) there exists at least one fusion channel \textit{c} in the fusion product $a \times a$ which possesses trivial spin factor, namely $\exists c \in a\times a$ such that $\theta_c=1$.

In general, anyon condensation breaks the topological order or symmetry  associated with a quantum group $\mathcal{H}$ down to a Hopf subalgebra $\mathcal{K}\subset \mathcal{H}$ \cite{BSS02,BSS02',BS08}. In other words, the irreducible representations of the quantum group $\mathcal{H}$ are exactly the charge sectors of the unbroken theory (or phase) and after condensation, the charge sectors of the broken theory carry the irreducible representations of $\mathcal{K}$. There are two general prescriptions under the anyon condensation: \textit{splitting} and \textit{identification}. The former refers to the fact that some irreps of the unbroken phase $\mathcal{H}$ will not be irreps of the broken phase $\mathcal{K}$ and splits into a number of irreps of $\mathcal{K}$ while the latter means some of the irreps of $\mathcal{H}$ will be equivalent representations under the action of $\mathcal{K}$ and hence be identified. These prescriptions can be encapsulated by the expression of the form:
\begin{equation}\label{Restriction of a}
a \rightarrow \sum_i n^i_a a_i,
\end{equation}
where $n^i_a\in \mathbb{Z}_{\geq 0}$ is the multiplicities of anyon species $a_i$. Note that $a$ and $a_i$ are an anyon (irrep) of $\mathcal{H}$ and $\mathcal{K}$ respectively.

The right-hand side of \eqref{Restriction of a} is called the restriction of $a$ and we say that the anyon $a$ restricts (or splits) to the anyon species $a_i$. On the other hand, we define the lift of an anyon $a_i$ of the broken phase as all the sectors $a's$ of the unbroken theory that contain $a_i$ in their restrictions. There are two assumptions to be imposed on \eqref{Restriction of a}. The first one is the set $C$ of condensed sectors must contain the new vacuum 1, that is $\forall c \in C$,
\begin{equation}
c\rightarrow 1+\sum_{i> 1}n^i_c c_i,
\end{equation}
where we have set $c_1=1$. In this sense, the condensate has been identified with the vacuum, although it may also split into other anyon species. Secondly, the fusion and splitting \eqref{Restriction of a} commute:
\begin{equation}
\left(\sum_k n^k_a a_k\right)\times \left(\sum_l n^l_b b_l\right)=\sum_{c,m} \tensor{N}{_{a,}_b^c} n^m_c c_m 
\end{equation}
The consequence is that the quantum dimensions are preserved under \eqref{Restriction of a}. 

Next, we consider the effect of confinement. In the broken phase, some of the charge sectors  will become confined and the others remain unconfined. If all the lifts of an anyon $a_i$ share the same spin factor, then $a_i$ becomes an unconfined anyon in the broken phase. Otherwise, it is a confined anyon in the condensed phase. Moreover, the set of unconfined sectors should be closed under fusion as well as has the vacuum as its element. Therefore, after considering confinement, we are left with a theory which only consists of unconfined sectors and we will say it is a phase after confinement.

In connection to our work, we start with an unbroken theory and it is described by an MFC, denoted as $\mathcal{B}$. After $\mathcal{B}$ undergone anyon condensation, we have a broken phase and denote it as $\mathcal{B}_{A}$, which stands for the category of right $A$-module. The theory which survives the confinement, on the other hand, is described by the category of dyslectic modules $\mathcal{B}^0_{A}$ which is a subcategory of $\mathcal{B}_A$.

\section{Classification}
\subsection{Definitions}
The monoidal products of the fusion categories $\mathcal{C}$'s are specified by fusion matrices $(N_i)_{jk}:={N_{i,j}}^k$ with $\mathbb N$-coefficients
\begin{align}
c_i\otimes c_j\cong\bigoplus_{k=1}^{\text{rank}(C)}{N_{i,j}}^k c_k.
\end{align}
Since the entries of fusion matrices are non-negative, we can apply the Perron-Frobenius theorem to get the largest eigenvalue, namely  Frobenius-Perron dimension of simple objects $c_i$, denoted as $\fp_{\mc C}(c_i)$. The Frobenius-Perron dimension of $\mathcal{C}$ is defined as 
\begin{align}\label{def_FPdimC}
\fp(\mc C):=\sum_{i=1}^{\text{rank}(\mc C)} (\fp_{\mc C}(c_i))^2,
\end{align}
namely the squared sum of Frobenius-Perron dimension of each simple objects. 

In an MFC, one can also define quantum dimension $d_i$ of $c_i$ by the quantum (or categorical)
trace
\begin{align*}
d_i:=\text{tr}(a_{c_i}),
\end{align*}
where $a:id_{\mc C}\cong (-)^{**}$ is a pivotal structure. Its multiplication rules are the same as the fusion rules of the simple objects
\begin{align}
d_id_j=\sum^{\text{rank}(\mc C)}_{k=1}{N_{i,j}}^kd_k.
\end{align}
The squared sum of $d_i$'s defines the categorical dimension
\begin{align}\label{D^C}
D^2(\mc C):=\sum^{\text{rank}(\mathcal{C})}_{i=1}d^2_i.
\end{align}
Hence, there are two $D(\mathcal{C})$, one positive and one negative, for each categorical dimension.

Additionally, a fusion category $\mathcal{C}$ can be equipped with a structure called braiding $c_{c_i,c_j}:c_i\otimes c_j\overset{\simeq}{\longrightarrow} c_j\otimes c_i$, making it a braided fusion category (BFC), with the braiding subject to the hexagon equations. (Here $c$'s denote generic objects in $\mathcal{C}$, which is a direct sum of simple objects $c_j$'s.) If the braiding is non-degenerate (and spherical), we have a modular fusion category (MFC). (For a physical introduction to the non-degeneracy of braiding, see \cite{K05}.)

The double braiding formula regarding two simple objects $b_i$ and $b_j$ is given by
\begin{equation}\label{Double Braiding}
c_{b_j,b_i}\cdot c_{b_i,b_j}\cong \sum_{k=1}^{\text{rank}(\mathcal{B})}\tensor{N}{_{i,j}^{k}}\frac{e^{2\pi i h_k}}{e^{2\pi i(h_i+h_j)}}id_k,
\end{equation}
where $id_k$ is the identity morphism at $b_k$. Taking the quantum trace of \eqref{Double Braiding} defines the unnormalized $S$-matrix
\begin{equation}
\tilde{S}_{i,j}:=\text{tr}(c_{b_j,b_i}\cdot c_{b_i,b_j})=\sum_{k=1}^{\text{rank}(\mathcal{B})}\tensor{N}{_{i,j}^{k}}\frac{e^{2\pi i h_k}}{e^{2\pi i(h_i+h_j)}}d_k,
\end{equation}
which, on the other hand, defines the normalized $S$-matrix (or just $S$-matrix)
\begin{equation}
S_{i,j}:=\frac{\tilde{S}_{i,j}}{D(\mathcal{B})}.
\end{equation} 
This defines an MFC in a precise mathematical sense: a spherical BFC (also called pre-modular FC) with non-degenerate $S$-matrix. Since we focus on MFC in this paper, the sets of quantum dimensions and conformal dimensions that lead to degenerate $S$-matrix should be ruled out.

The $S$-matrix obeys one of the important relations
\begin{equation}\label{Square of S=C}
S^2=C,
\end{equation}
where the matrix $C$ is the charge conjugation matrix, defined by 
\begin{equation}
C_{i,j}=\delta_{i,j}\quad (b^{*}_i\cong b_j),
\end{equation}
where $b^*_i$ is the dual of $b_i$. The $S$-matrix also satisfies 
\begin{equation}
\widetilde{S}_{i,j^{*}}=\left(\widetilde{S}_{i,j}\right)^*,
\end{equation}
where the RHS stands for complex conjugation. In general, we use \eqref{Square of S=C} to find the permissible conformal dimensions. However, we have a more restrictive situation if the MFC is self-dual, meaning all its simple objects are self-dual. Then the elements of $S$-matrix should be real
\begin{equation}
\forall b_j\in \mathcal{B},\ S_{i,j}\in \mathbb{R} \quad\text{if}\quad b^*_i\cong b_i.
\end{equation} 
One could use this to find the set of $h_j$'s instead.

Next we review definitions on algebras (in an MFC). An algebra in a fusion category $\mathcal{C}$ is a triple $(A,\mu,u)$ of an object $A\in \mathcal{C}$, multiplication morphism $\mu:A\otimes A\rightarrow A$, and unit morphism $u:1\rightarrow A$ obeying associativity and unit axioms. A category of right $A$-modules consists of pairs $(m,p)$ where $m\in\mathcal{C}$ and $p:m\otimes A \rightarrow m$ subject to consistency conditions. An algebra is called separable if $\mathcal{C}_A$ is semisimple. 

An algebra $A\in \mathcal{B}$ in a BFC is called commutative if 
\begin{equation}\label{Commutatitive(1)}
\mu\cdot c_{A,A}=\mu.
\end{equation}
A commutative separable algebra is called \'etale. An algebra $A\in\mc C$ is called connected if $\text{dim}_{\mathbb{C}}\mathcal{C}(1,A)=1$. A connected \'etale algebra $A\in \mathcal{B}$ is called Lagrangian if $\Bigl(\fp_{\mc B}(A)\Bigr)^2=\fp(\mathcal{B})$. A BFC without non-trivial connected \'etale algebra is called completely anisotropic. The category of right $A$-modules $\mathcal{B}_A$ contains an important subcategory $ \mathcal{B}^0_A$, which consists of dyslectic (or local) modules \cite{P95} $(m,p)\in \mathcal{B}_A$ obeying
\[ p\cdot c_{A,m}\cdot c_{m,A}=p.\]

The category $\mathcal{C}_A$ of right $A$-modules is a left $\mathcal{C}$-module category \cite{O01}. A left $\mathcal{C}$-module category is a quadruple $(\mathcal{M},\triangleright,m,l)$ of a category $\mathcal{M}$, an action bifunctor $\triangleright: \mc C\times \mc M\rightarrow \mc M$, natural isomorphisms $m_{-,-,-}:(-\otimes -)\triangleright - \cong - \triangleright (- \triangleright -)$ and $l:1\triangleright \mc M\simeq \mc M$ called module associativity constraint and unit constraint respectively. Both of them have to satisfy the associativity and unit axioms.

\subsection{Method}
It turns out that the modularity of $\mathcal{B}$ largely reduces the computational work required in classifying the \'etale algebras. Here we aim to present the methods used in this paper in a self-contained manner. For the scheme which applies to the more general case of braided fusion categories, see \cite{KK23preMFC,KK23rank5}. We review here the three necessary conditions that a BFC $\mathcal{B}$ has to satisfy (the symbol is the same for a MFC, but it should be clear in the context); there should exist a fusion category $\mc C$ such that
\begin{align}\label{3 necessary conditions for BFC}
(\text{i})\, &1\leq \text{rank}(\mathcal{C})\leq \lfloor \fp (\mathcal{B})\rfloor,\nonumber\\ 
(\text{ii})\,&1\leq \fp(\mathcal{C})\leq \fp(\mathcal{B}),\nonumber\\ 
(\text{iii})\,&\fp_{\mathcal{B}}(A)=\frac{\fp(B)}{\fp(\mathcal{C})}.
\end{align} 
When $A$ is connected \'etale, we get a fusion category $\mc B_A$ with Frobenius-Perron dimension
\begin{equation}
    \fp(\mc B_A)=\frac{\fp(\mc B)}{\fp_\mc B(A)}.\label{FPdimBA}
\end{equation}
The key point is if further $\mathcal{B}$ is  modular and $A\in \mathcal{B}$ is a connected \'etale algebra, then $\mathcal{B}^0_{A}$ is modular and obeys \cite{P95,KO01,EGNO15}
\begin{equation}\label{FPdim B^0_A}
\fp(\mathcal{B}^0_{A})= \frac{\fp(\mc B)}{\Bigl(\fp_{\mc B}(A)\Bigr)^2}.
\end{equation}
We also have \cite{ENO02,EGNO15}
\begin{equation}\label{FPdim larger than 1}
\forall c\in \mathcal{C},\quad \fp_{\mc C}(c)\geq 1.
\end{equation}
This can be proved by using the Perron-Frobenius theorem and the fact $(N_i)_{j,k}\in \mathbb{N}$ (see Appendix A of \cite{NHKSB16}).

Using \eqref{FPdim B^0_A} and \eqref{FPdim larger than 1}, we may reach
\begin{equation} 
    1\le\left(\fp_{\mc B}(A)\right)^2\le\fp(\mc B).\label{FPdimA2bound}
\end{equation}
The squaring of Frobenius-Perron dimension of $A$ in \eqref{FPdimA2bound} poses a greater constraint than the general case mentioned above, in particular the condition (ii) in \eqref{3 necessary conditions for BFC} and consequently less candidates for $A$ will be found. By definition, we can take an ansatz for $A$ in the most general form
\begin{equation}\label{Ansatz for A}
A\cong 1\oplus \bigoplus_{j\neq 1} n_j b_j,
\end{equation}
where $n_j\in \mathbb{N}$ and we have set $n_1=1$ as well as $b_1\cong 1$, implying $A$ contains the unique identity. It has 
\begin{equation}
\fp_{\mathcal{B}}(A)=1+\sum_{j\neq 1}n_j\fp_{\mc B}(b_j).
\end{equation}

Then we insert \eqref{Ansatz for A} into \eqref{FPdimA2bound} to find out the sets of $n_j\text{'}s$ which specifies the list of candidates and check if each candidate satisfies the conditions of connected \'etale algebra. To achieve this, three conditions are to be satisfied for $A$: connectedness, separability and commutativity. The connectedness condition is evident from \eqref{Ansatz for A} whereas the separability is guaranteed as long as $\mathcal{B}_A$ is a fusion category. On the other hand, the commutativity condition
\begin{equation}\label{Commutativity(2)}
\mu\cdot c_{A,A}=\mu,
\end{equation}
is non-trivial to prove. However, we may simply check the necessary condition at the first place, given by
\begin{equation}
\mu\cdot c_{A,A}\cdot c_{A,A}= \mu.
\end{equation}
Using the formula \eqref{Double Braiding}, it is sufficient to check instead
\begin{equation}\label{Necessary Conditions For etale}
c_{A,A}\cdot c_{A,A}\cong \sum_{i,j=1}^{\text{rank}(\mathcal{B})}n_i n_j (\iota_i\otimes\iota_j)\cdot c_{b_j,b_i}\cdot c_{b_i,b_j}\cdot(p_i\otimes p_j).
\end{equation}
In many cases, the necessary condition \eqref{Necessary Conditions For etale} is enough for ruling out many candidates. Otherwise, we may turn to \eqref{Commutatitive(1)} or \eqref{Commutativity(2)} for final settlement. 

For each non-trivial \'etale algebra $A$ in $\mathcal{B}$, it is of interest to find out the correct $\mathcal{B}_A$. As mentioned above, the anyon condensation provides us a means to do so. For an illustration of this method, see section \ref{so52} for $so(5)_2$ as an example, the only multiplicity-free rank six MFC with non-trivial connected \'etale algebra given by $1\oplus X$. There it is found that the phase after confinement is $\vecG^1_{\mathbb{Z}/5\mathbb{Z}}$. Here we remark that anyon condensation is only one of the ways in identifying $\mathcal{B}_A$.

As a demonstration of our method, we will present the toric code as an example, where $\mc B\simeq \tc$. The four simple objects $\{1,X,Y,Z\}$  have the following fusion rules:
\begin{table}[H]
\begin{center}
\begin{tabular}{c|c|c|c|c}
    $\otimes$&$1$&$X$&$Y$&$Z$\\\hline
    $1$&$1$&$X$&$Y$&$Z$\\\hline
    $X$&&$1$&$Z$&$Y$\\\hline
    $Y$&&&$1$&$X$\\\hline
    $Z$&&&&$1$
\end{tabular}.
\end{center}
\end{table}

Therefore, the Frobenius-Perron dimensions are \[ \hspace{-30pt}\fp_{\mcal B}(1)=\fp_{\mcal B}(X)=\fp_{\mcal B}(Y)=\fp_{\mcal B}(Z)=1. \]
The quantum dimension of each simple object follows the same fusion rules, leading to the conditions: $d^2_X=d^2_Y=d^2_Z=1,\, d_Xd_Y=d_Z,\, d_Yd_Z=d_X,\, d_Zd_X=d_Y$. The solutions for $(d_X,d_Y,d_Z)$ are:
\begin{align*}
    (d_X,d_Y,d_Z)=(1,1,1),(1,-1,-1),(-1,1,-1),(-1,-1,1).
\end{align*}
From (\ref{def_FPdimC}), we find that
\begin{align*}
    \fp(\mc B)=4.
\end{align*}
Thus, according to (\ref{FPdimA2bound}), we can write the bound:
\begin{align*}
    1\leq\left(\fp_{\mcal B}(A)\right)^2\leq 4.
\end{align*}
The general ansatz for the connected algebra is given by
\[ A\cong1\oplus n_XX\oplus n_YY\oplus n_ZZ, \]
which implies
\[ \fp_{\mc B}(A)=1+n_X+n_Y+n_Z. \]
The natural numbers $(n_X,n_Y,n_Z)$ can take the values:
\begin{align*}
    (n_X,n_Y,n_Z)=(0,0,0),(1,0,0),(0,1,0),(0,0,1),\\
\end{align*}
corresponding to
\begin{align*}
    A\cong 1,\, 1\oplus X,\, 1\oplus Y,\, 1\oplus Z.\\
\end{align*}
For the conformal dimensions, the possible values of $(h_X,h_Y,h_Z)$ are:
\begin{align*}
    (h_X,h_Y,h_Z)=
    \begin{cases}
        (\frac{1}{2},0,0),(\frac{1}{2},\frac{1}{2},\frac{1}{2}) & (d_X,d_Y,d_Z)=(1,1,1) \\
        (\frac{1}{2},0,0),(0,\frac{1}{2},0),(\frac{1}{2},\frac{1}{2},\frac{1}{2}) & (d_X,d_Y,d_Z)=(1,-1,-1)
    \end{cases}
\end{align*}
In order for $A$ to become connected \'etale algebra, the following conditions must be satisfied:

1.\textbf{Algebra} \\
Based on the results in \cite{KK23preMFC}, $1\oplus X$, as a $\mbb Z/2\mbb Z$ object, is an algebra in $\vecG_{\mbb Z/2\mbb Z}$, which is a subcategory of $\mc B$. Consequently, by the lemma in \cite{rank9}, $1\oplus X$ is also an algebra in $\mc B$. This applies similarly to $1\oplus Y$ and $1\oplus Z$, since both $Y$ and $Z$ are $\mbb Z/2\mbb Z$ objects.

2.\textbf{Connected} \\
The coefficient of the object $1$ in the algebra must be 1. A generic \'etale algebra is a direct sum of connected \'etale algebras, so it suffices to study the connected \'etale algebras.

3.\textbf{Commutative} \\
According to (\ref{Commutativity(2)}), we need $c_{A,A}\cong id_{A\otimes A}$, and this has been proved in \cite{KK23preMFC}.

4.\textbf{Separable} \\
Since $\mc B_A$ is a fusion category (which we will prove it later), it follows that $\mc B_A$ is semisimple. Hence, $A$ is separable.

From the result in \cite{KK23preMFC}, given a simple object $X$, in order for $1\oplus X$ to be connected \'etale algebra, we need $(d_X,h_X)=(1,0)$. Therefore, we obtain the following candidate algebras: 
\begin{align*}
     A\cong
     \begin{cases}
        1 & \text{all MFCs}\\
        1\oplus Y,\;1\oplus Z& (d_X,d_Y,d_Z)=(1,1,1),\; (h_X,h_Y,h_Z)=(\frac{1}{2},0,0) \\
        1\oplus X& (d_X,d_Y,d_Z)=(1,-1,-1),\; (h_X,h_Y,h_Z)=(0,\frac{1}{2},0)
    \end{cases}.
\end{align*}
In order to identify $\mc B_A$, we have to compute new fusion rule in $\mc B_A$. It turns out the free module functor is useful. Here, the free module functor $F_A:\mc B\to\mc B_A$ is defined by
\begin{align*}
    F_A(b_i):=b_i\otimes A,
\end{align*}
where $b_i$ is a simple object in $\mc B$. The functor satisfies the following property:
\begin{align*}
    F_A(b_i)\otimes_A F_A(b_j)\cong F_A(b_i\otimes b_j).
\end{align*}
Now, for the first case where $A\cong 1\oplus Y$, we have:
\begin{gather*}
    F_A(1)\cong1\oplus Y,\quad F_A(X)\cong X\oplus Z,\quad F_A(Y)\cong1\oplus Y,\quad F_A(Z)\cong X\oplus Z.
\end{gather*}
The only distinct simple objects are
\[ m_1\cong 1\oplus Y,\quad m_2\cong X\oplus Z. \]
The fusion rules of them can be computed with the free module functor. For example, we have
\begin{align*}
    m_2\otimes_Am_2&\equiv F_A(X)\otimes_AF_A(X)\\
    &\cong F_A(X\otimes X)\\
    &\cong F_A(1)\equiv m_1.
\end{align*}
Similarly, one can compute all fusion rules
\begin{table}[H]
\begin{center}
\begin{tabular}{c|c|c}
    $\otimes_A$&$m_1$&$m_2$\\\hline
    $m_1$&$m_1$&$m_2$\\\hline
    $m_2$&$m_2$&$m_1$
\end{tabular}.
\end{center}
\end{table}
\hspace{-17pt}This shows that $\mc B_A\simeq \vecG_{\mbb Z/2\mbb Z}$. Since $\vecG_{\mbb Z/2\mbb Z}$ is semisimple, $A$ is separable.
From the fusion table, we determine that $\fp(\mc B_A)=2$. This agrees with the formula (\ref{FPdimBA}). Using (\ref{FPdim B^0_A}), we also find $\fp(\mc B_A^0)=1$, leading to the identification:
\begin{align*}
    \mc B_A^0\simeq\vect=\{m_1\}.
\end{align*}

For the second case, where $A\cong 1\oplus X$, we have:
\begin{align*}
    F_A(1)\cong1\oplus X\cong F_A(X),\quad F_A(Y)\cong Y\oplus Z\cong F_A(Z).
\end{align*}
We label
\[
m_1 \cong 1\oplus X,\quad m_2\cong Y\oplus Z ,
\]
which satisfy the following fusion table:
\begin{table}[H]
\begin{center}
\begin{tabular}{c|c|c}
    $\otimes_A$&$m_1$&$m_2$\\\hline
    $m_1$&$m_1$&$m_2$\\\hline
    $m_2$&$m_2$&$m_1$
\end{tabular}.
\end{center}
\end{table}
\hspace{-17pt}This yields the same result, $\mc B_A\simeq\vecG_{\mbb Z/2\mbb Z}$. We also identify $\mc B_A^0\simeq\vect=\{m_1\}$.

\subsection{Classification}
\subsubsection{$\mc B\simeq\vecG
_{\mbb Z/6\mbb Z}^\alpha$}
The fusion table is given by
\begin{table}[H]
\begin{center}
\begin{tabular}{c|c|c|c|c|c|c}
    $\otimes$&$1$&$X$&$Y$&$Z$&$V$&$W$\\\hline
    $1$&$1$&$X$&$Y$&$Z$&$V$&$W$\\\hline
    $X$&&$1$&$W$&$V$&$Z$&$Y$\\\hline
    $Y$&&&$V$&$1$&$X$&$Z$\\\hline
    $Z$&&&&$W$&$Y$&$X$\\\hline
    $V$&&&&&$W$&$1$\\\hline
    $W$&&&&&&$V$
\end{tabular}.
\end{center}
\end{table}
\hspace{-17pt}The Frobenius-Perron dimensions are all the same
\[ \hspace{-30pt}\fp_{\mcal B}(1)=\fp_{\mcal B}(X)=\fp_{\mcal B}(Y)=\fp_{\mcal B}(Z)=\fp_{\mcal B}(V)=\fp_{\mcal B}(W)=1, \]
and
\[ \fp(\mc B)=6. \]
The quantum dimensions $d_j$'s obey the equations
$d_X^2=1,d_Xd_Y=d_W,d_Xd_Z=d_V,d_Xd_V=d_Z,d_Xd_W=d_Y,d_Y^2=d_V,d_Yd_Z=1,
d_Yd_V=d_X,d_Yd_W=d_Z,d_Z^2=d_W,d_Zd_V=d_Y,d_Zd_W=d_X,
d_V^2=d_W,d_Vd_W=1,d_W^2=d_V$. There are two solutions
\[ (d_X,d_Y,d_Z,d_V,d_W)=(-1,-1,-1,1,1),(1,1,1,1,1). \]
The second solution gives unitary MFCs. They have the same categorical dimension
\[ D^2(\mc B)=6. \]
Regardless of quantum dimensions, there are 4 conformal dimensions
$$(h_X,h_Y,h_Z,h_V,h_W)=(\frac{1}{4},\frac{7}{12},\frac{7}{12},\frac{1}{3},\frac{1}{3}),(\frac{1}{4},\frac{11}{12},\frac{11}{12},\frac{2}{3},\frac{2}{3}),(\frac{3}{4},\frac{1}{12},\frac{1}{12},\frac{1}{3},\frac{1}{3}),(\frac{3}{4},\frac{5}{12},\frac{5}{12},\frac{2}{3},\frac{2}{3}).$$
Therefore, there are
\[ 2(\text{quantum dimensions})\times4(\text{conformal dimensions})\times2(\text{categorical dimensions})=16 \]
MFCs, among which those eight with the last quantum dimensions give unitary MFCs. We classify connected étale algebras in all 16 MFCs simultaneously.

There are six sets of natural numbers that obey (\ref{FPdimA2bound})
\begin{align*}
    (n_X,n_Y,n_Z,n_V,n_W)=&(0,0,0,0,0),(1,0,0,0,0),(0,1,0,0,0),\\
    &(0,0,1,0,0),(0,0,0,1,0),(0,0,0,0,1).
\end{align*}
The first solution is nothing but the trivial connected étale algebra $A\cong1$ giving $\mc B_A^0\simeq\mc B_A\simeq\mc B$. All the other solutions do not lead to commutative algebra since they include object(s) with nontrivial conformal dimensions and hence the necessary condition is not fulfilled.

Therefore, we have
\begin{table}[H]
\begin{center}
\begin{tabular}{c|c|c|c}
    Connected étale algebra $A$&$\mc B_A$&$\rank(\mc B_A)$&Lagrangian?\\\hline
    $1$&$\mc B$&$6$&No
\end{tabular}.
\end{center}
\caption{Connected étale algebras in rank six MFC $\mcal B\simeq \vecG_{\mbb Z/6\mbb Z}^\alpha$}\label{rank6Z6results}
\end{table}
\hspace{-17pt}

\subsubsection{$\mc B\simeq\vecG_{\mbb Z/2\mbb Z}^{-1}\boxtimes\ising$}
The MFC is given by a Deligne tensor product \cite{D90}, whose fusion table is given by
\begin{table}[H]
\begin{center}
\begin{tabular}{c|c|c|c|c|c|c}
    $\otimes$&$1$&$X$&$Y$&$Z$&$V$&$W$\\\hline
    $1$&$1$&$X$&$Y$&$Z$&$V$&$W$\\\hline
    $X$&&$1$&$Z$&$Y$&$W$&$V$\\\hline
    $Y$&&&$1$&$X$&$W$&$V$\\\hline
    $Z$&&&&$1$&$V$&$W$\\\hline
    $V$&&&&&$1\oplus Z$&$X\oplus Y$\\\hline
    $W$&&&&&&$1\oplus Z$

\end{tabular}.
\end{center}
\end{table}
The Frobenius-Perron dimensions are
\[ \fp_{\mc B}(1)=\fp_{\mc B}(X)=\fp_{\mc B}(Y)=\fp_{\mc B}(Z)=1,\fp_{\mc B}(V)=\fp_{\mc B}(W)=\sqrt2 \]
and
\[ \fp(\mc B)=8. \]
Then, the solutions of quantum dimension are 
\begin{align*}
&(d_X,d_Y,d_Z,d_V,d_W)&\\
&=
(1,1,1,\sqrt{2},\sqrt{2}),
(-1,-1,1,\sqrt{2},-\sqrt{2}),
(-1,-1,1,-\sqrt{2},\sqrt{2}),
(1,1,1,-\sqrt{2},-\sqrt{2}).&
\end{align*}
Therefore, the categorical dimension is
\[ D^2(\mc B)=8.\]
They have the following conformal dimensions. For the first and fourth quantum dimensions, the two conformal dimensions related by permutations $(XY)$ or $(VW)$ lead to the same MFC. Thus different MFCs are given by eight conformal dimensions
\begin{align*}
(h_X,h_Y,h_Z,h_V,h_W)&=(\frac14,\frac34,\frac12,\frac1{16},\frac5{16}),(\frac14,\frac34,\frac12,\frac1{16},\frac{13}{16}),(\frac14,\frac34,\frac12,\frac3{16},\frac7{16}),(\frac14,\frac34,\frac12,\frac3{16},\frac{15}{16}),\\
&~~~~(\frac14,\frac34,\frac12,\frac5{16},\frac9{16}),(\frac14,\frac34,\frac12,\frac7{16},\frac{11}{16}),(\frac14,\frac34,\frac12,\frac9{16},\frac{13}{16}),(\frac14,\frac34,\frac12,\frac{11}{16},\frac{15}{16})\quad(\rm{mod\ 1}).
\end{align*}
On the other hand, for the second and third quantum dimensions, permutations $(VW)$ give different MFCs. Thus, different MFCs are given by 16 conformal dimensions
\begin{align*}
(h_X,h_Y,h_Z,h_V,h_W)&=(\frac14,\frac34,\frac12,\frac1{16},\frac5{16}),(\frac14,\frac34,\frac12,\frac1{16},\frac{13}{16}),(\frac14,\frac34,\frac12,\frac3{16},\frac7{16}),(\frac14,\frac34,\frac12,\frac3{16},\frac{15}{16}),\\
&~~~~(\frac14,\frac34,\frac12,\frac5{16},\frac9{16}),(\frac14,\frac34,\frac12,\frac7{16},\frac{11}{16}),(\frac14,\frac34,\frac12,\frac9{16},\frac{13}{16}),(\frac14,\frac34,\frac12,\frac{11}{16},\frac{15}{16}),\\
&~~~~(\frac14,\frac34,\frac12,\frac5{16},\frac1{16}),(\frac14,\frac34,\frac12,\frac{13}{16},\frac1{16}),(\frac14,\frac34,\frac12,\frac7{16},\frac3{16}),(\frac14,\frac34,\frac12,\frac{15}{16},\frac3{16}),\\
&~~~~(\frac14,\frac34,\frac12,\frac9{16},\frac5{16}),(\frac14,\frac34,\frac12,\frac{11}{16},\frac7{16}),(\frac14,\frac34,\frac12,\frac{13}{16},\frac9{16}),(\frac14,\frac34,\frac12,\frac{15}{16},\frac{11}{16})\quad(\rm{mod\ 1}).
\end{align*}
\newline

Therefore, with the two signs of categorical dimensions, there are
\[ 16+16+32+32=96 \]
MFCs, among which those 16 with the first quantum dimensions are unitary. We classify connected étale algebras in all 96 MFCs simultaneously.

The set of natural numbers for the ansatz obeying (\ref{FPdimA2bound}) are 
\begin{align*}
(n_X,n_Y,n_Z,n_V,n_W)=
(0,0,0,0,0),(1,0,0,0,0),(0,1,0,0,0),\\
(0,0,1,0,0),(0,0,0,1,0),(0,0,0,0,1).
\end{align*}
The only candidate of connected étale algebra is $A\cong1$ giving $\mc B_A^0\cong\mc B_A\cong\mc B.$\\
We conclude
\begin{table}[H]
\begin{center}
\begin{tabular}{c|c|c|c}
    Connected étale algebra $A$&$\mc B_A$&$\rank(\mc B_A)$&Lagrangian?\\\hline
    $1$&$\mc B$&$6$&No
\end{tabular}.
\end{center}
\caption{Connected étale algebras in rank six MFC $\mcal B\simeq\vecG_{\mbb Z/2\mbb Z}^{-1}\boxtimes\ising$}\label{rank6BFR601results}
\end{table}

\subsubsection{$\mc B\simeq su(3)_2\simeq\fib\boxtimes\vecG_{\mbb Z/3\mbb Z}^1$}
The MFC has six simple objects $\{1,X,Y,Z,V,W\}$ obeying monoidal products
\begin{table}[H]
\begin{center}
\begin{tabular}{c|c|c|c|c|c|c}
    $\otimes$&$1$&$X$&$Y$&$Z$&$V$&$W$\\\hline
    $1$&$1$&$X$&$Y$&$Z$&$V$&$W$\\\hline
    $X$&&$Y$&$1$&$W$&$Z$&$V$\\\hline
    $Y$&&&$X$&$V$&$W$&$Z$\\\hline
    $Z$&&&&$1\oplus Z$&$Y\oplus V$&$X\oplus W$\\\hline
    $V$&&&&&$X\oplus W$&$1\oplus Z$\\\hline
    $W$&&&&&&$Y\oplus V$
\end{tabular}.
\end{center}
\end{table}
\hspace{-17pt}Thus, they have
\[ \hspace{-30pt}\fp_{\mcal B}(1)=\fp_{\mcal B}(X)=\fp_{\mcal B}(Y)=1,\quad \fp_{\mcal B}(Z)=\fp_{\mcal B}(V)=\fp_{\mcal B}(W)=\zeta:=\frac{1+\sqrt{5}}{2}, \]
and
\[ \fp(\mc B)=\frac{15+3\sqrt5}2\approx10.9. \]
There are two solutions of quantum dimension
\[ (d_X,d_Y,d_Z,d_V,d_W)=(1,1,-\zeta^{-1},-\zeta^{-1},-\zeta^{-1}),(1,1,\zeta,\zeta,\zeta).  \]
The last solution yields unitary MFCs. They have categorical dimensions
\[ D^2(\mc B)=\frac{15-3\sqrt{5}}{2}(\approx 4.1),\frac{15+3\sqrt{5}}{2}(\approx 10.9). \]
The conformal dimensions are
\[ (h_X,h_Y,h_Z,h_V,h_W)=\begin{cases}(\frac13,\frac13,\frac15,\frac{8}{15},\frac{8}{15}),(\frac13,\frac13,\frac45,\frac{2}{15},\frac{2}{15}),(\frac23,\frac23,\frac15,\frac{13}{15},\frac{13}{15}),(\frac23,\frac23,\frac45,\frac{7}{15},\frac{7}{15})&(d_Z=-\zeta^{-1}),\\(\frac13,\frac13,\frac25,\frac{11}{15},\frac{11}{15}),(\frac13,\frac13,\frac35,\frac{14}{15},\frac{14}{15}),(\frac23,\frac23,\frac25,\frac{1}{15},\frac{1}{15}),(\frac23,\frac23,\frac35,\frac{4}{15},\frac{4}{15})&(d_Z=\zeta).\end{cases}\]
Therefore, there are 
\[ 2(\text{quantum dimensions})\times4(\text{conformal dimensions})\times2(\text{categorical dimensions})=16 \]
MFCs, among which those eight with the last quantum dimensions give unitary MFCs. We classify connected étale algebras in all 16 MFCs simultaneously.

The natural numbers obeying (\ref{FPdimA2bound}) are
\begin{align*}
    (n_X,n_Y,n_Z,n_V,n_W)=&(0,0,0,0,0),(1,0,0,0,0),(0,1,0,0,0),(0,0,1,0,0)\\
    &(0,0,0,1,0),(0,0,0,0,1),(1,1,0,0,0),(2,0,0,0,0)\\
    &(0,2,0,0,0).
\end{align*}
The first solution is nothing but the trivial connected étale algebra $A\cong1$ giving $\mc B_A^0\simeq\mc B_A\simeq\mc B$. The second to sixth candidates fail to give commutative algebra since they all possess nontrivial conformal dimensions. Then the non-commutativity of the second candidate is sufficient to guarantee the seventh and eighth candidates fail to pass the necessary condition as well. Similar reasoning applies to the third and the final candidate. In other words, $\mc B\simeq su(3)_2\simeq\fib\boxtimes\vecG_{\mbb Z/3\mbb Z}^1$ is completely anisotropic.

Summarizing, we have
\begin{table}[H]
\begin{center}
\begin{tabular}{c|c|c|c}
    Connected étale algebra $A$&$\mc B_A$&$\rank(\mc B_A)$&Lagrangian?\\\hline
    $1$&$\mc B$&$6$&No
\end{tabular}.
\end{center}
\caption{Connected étale algebras in rank six MFC $\mc B \simeq\fib\boxtimes\vecG_{\mbb Z/3\mbb Z}^1$}\label{rank6su32results}
\end{table}
\hspace{-17pt}

\subsubsection{$\mc B\simeq\text{TriCritIsing}$}
The MFC has six simple objects $\{1,X,Y,Z,V,W\}$ obeying monoidal products
\begin{table}[H]
\begin{center}
\begin{tabular}{c|c|c|c|c|c|c}
    $\otimes$&$1$&$X$&$Y$&$Z$&$V$&$W$\\\hline
    $1$&$1$&$X$&$Y$&$Z$&$V$&$W$\\\hline
    $X$&&$1$&$Y$&$V$&$Z$&$W$\\\hline
    $Y$&&&$1\oplus X$&$W$&$W$&$Z\oplus V$\\\hline
    $Z$&&&&$1\oplus V$&$X\oplus Z$&$Y\oplus W$\\\hline
    $V$&&&&&$1\oplus V$&$Y\oplus W$\\\hline
    $W$&&&&&&$1\oplus X\oplus Z\oplus V$

\end{tabular}.
\end{center}
\end{table}
Thus, they have

\begin{align*}
\fp_{\mc B}(1)=1=\fp_{\mc B}(X),&\quad\fp_{\mc B}(Y)=\sqrt{2},\\\fp_{\mc B}(Z)=\frac{1+\sqrt5}2=\fp_{\mc B}(V),&\quad\fp_{\mc B}(W)=\sqrt{3+\sqrt{5}}.
\end{align*}

and
\[ \fp(\mc B)=10+2\sqrt{5}(\approx14.5). \]

The quantum dimensions are given by 
\begin{align*}
&(d_X,d_Y,d_Z,d_V,d_W)&\\
=&
(1,\sqrt{2},\frac{1-\sqrt{5}}2,\frac{1-\sqrt{5}}2,-\sqrt{3-\sqrt{5}}),
(1,-\sqrt{2},\frac{1-\sqrt{5}}2,\frac{1-\sqrt{5}}2,
\sqrt{3-\sqrt{5}}),\\
&(1,\sqrt{2},\frac{1+\sqrt{5}}2,\frac{1+\sqrt{5}}2,\sqrt{3+\sqrt{5}}),
(1,-\sqrt{2},\frac{1+\sqrt{5}}2,\frac{1+\sqrt{5}}2,-\sqrt{3+\sqrt{5}})
\end{align*}
with categorical dimensions
\[ D^2(\mc B)=10-2\sqrt5(\approx5.5),\quad10+2\sqrt5,\]
respectively for each pair. Meanwhile they have 16 conformal dimensions
\begin{align*}
(h_X,h_Y,h_Z,h_V,h_W)=
\begin{cases}
(\frac12,\frac1{16},\frac9{10},\frac25,\frac{37}{80}),(\frac12,\frac1{16},\frac1{10},\frac35,\frac{53}{80}),(\frac12,\frac3{16},\frac9{10},\frac25,\frac{47}{80}),(\frac12,\frac3{16},\frac1{10},\frac35,\frac{63}{80}),\\
(\frac12,\frac5{16},\frac9{10},\frac25,\frac{57}{80}),(\frac12,\frac5{16},\frac1{10},\frac35,\frac{73}{80}),(\frac12,\frac7{16},\frac9{10},\frac25,\frac{67}{80}),(\frac12,\frac7{16},\frac1{10},\frac35,\frac{3}{80}),\\
(\frac12,\frac{9}{16},\frac9{10},\frac25,\frac{77}{80}),(\frac12,\frac{9}{16},\frac1{10},\frac35,\frac{13}{80}),(\frac12,\frac{11}{16},\frac9{10},\frac25,\frac{7}{80}),(\frac12,\frac{11}{16},\frac1{10},\frac35,\frac{23}{80}),\\
(\frac12,\frac{13}{16},\frac9{10},\frac25,\frac{17}{80}),(\frac12,\frac{13}{16},\frac1{10},\frac35,\frac{33}{80}),(\frac12,\frac{15}{16},\frac9{10},\frac25,\frac{27}{80}),
(\frac12,\frac{15}{16},\frac1{10},\frac35,\frac{43}{80}),
\end{cases}
(\rm{mod\ 1})
\end{align*}
for $d_V=\frac{1+\sqrt5}2$, and
\begin{align*}
(h_X,h_Y,h_Z,h_V,h_W)=
\begin{cases}
(\frac12,\frac1{16},\frac7{10},\frac15,\frac{21}{80}),(\frac12,\frac1{16},\frac3{10},\frac45,\frac{69}{80}),(\frac12,\frac3{16},\frac7{10},\frac15,\frac{31}{80}),(\frac12,\frac3{16},\frac3{10},\frac45,\frac{79}{80}),\\
(\frac12,\frac5{16},\frac7{10},\frac15,\frac{41}{80}),(\frac12,\frac5{16},\frac3{10},\frac45,\frac{9}{80}),(\frac12,\frac7{16},\frac7{10},\frac15,\frac{51}{80}),(\frac12,\frac7{16},\frac3{10},\frac45,\frac{19}{80}),\\
(\frac12,\frac{9}{16},\frac7{10},\frac15,\frac{61}{80}),(\frac12,\frac{9}{16},\frac3{10},\frac45,\frac{29}{80}),(\frac12,\frac{11}{16},\frac7{10},\frac15,\frac{71}{80}),(\frac12,\frac{11}{16},\frac3{10},\frac45,\frac{39}{80}),\\
(\frac12,\frac{13}{16},\frac7{10},\frac15,\frac{1}{80}),(\frac12,\frac{13}{16},\frac3{10},\frac45,\frac{49}{80}),(\frac12,\frac{15}{16},\frac7{10},\frac15,\frac{11}{80}),
(\frac12,\frac{15}{16},\frac3{10},\frac45,\frac{59}{80}),
\end{cases}
(\rm{mod\ 1})
\end{align*}
for $d_V=\frac{1-\sqrt{5}}2$.\newline

Thus, there are
\[ 4(\text{quantum dimensions})\times16(\text{conformal dimensions})\times2(\text{categorical dimensions})=128 \]
MFCs, among which those 32 with the third quantum dimensions are unitary.
The general candidate of connected étale algebra can be written as
\[ A\cong1\oplus n_XX\oplus n_YY\oplus n_ZZ\oplus n_VV\oplus n_WW. \]
with $n_j\in\mathbb{N}$. It has 
\[
\fp_{\mc B}(A)=1+n_X+\sqrt{2}n_Y+\dfrac{1+\sqrt{5}}{2}n_Z+\dfrac{1+\sqrt{5}}{2}n_V+\sqrt{3+\sqrt{5}}n_W.
\]
In order to obey the necessary conditions of étale algebra (\ref{FPdimA2bound}), the only possible sets are given by
\begin{align*}
(n_X,n_Y,n_Z,n_V,n_W)=
&(0,0,0,0,0),(1,0,0,0,0),(0,1,0,0,0),\\
&(0,0,1,0,0),(0,0,0,1,0),(0,0,0,0,1),\\
&(1,1,0,0,0),(1,0,1,0,0),(1,0,0,1,0),\\
&(2,0,0,0,0).
\end{align*}
The only candidate of connected étale algebra is $A\cong1$ giving $\mc B_A^0\simeq\mc B_A\simeq\mc B.$

We conclude
\begin{table}[H]
\begin{center}
\begin{tabular}{c|c|c|c}
    Connected étale algebra $A$&$\mc B_A$&$\rank(\mc B_A)$&Lagrangian?\\\hline
    $1$&$\mc B$&$6$&No
\end{tabular}.
\end{center}
\caption{Connected étale algebras in rank six MFC $\mcal B\simeq\text{TriCritIsing}$}\label{rank6tricritisingresults}
\end{table}

\subsubsection{$\mc B\simeq su(2)_5\simeq\vecG_{\mbb Z/2\mbb Z}^{-1}\boxtimes psu(2)_5$}
The MFC has six simple objects $\{1,X,Y,Z,V,W\}$ obeying monoidal products
\begin{table}[H]
\begin{center}
\begin{tabular}{c|c|c|c|c|c|c}
    $\otimes$&$1$&$X$&$Y$&$Z$&$V$&$W$\\\hline
    $1$&$1$&$X$&$Y$&$Z$&$V$&$W$\\\hline
    $X$&&$1$&$Z$&$Y$&$W$&$V$\\\hline
    $Y$&&&$1\oplus V$&$X\oplus W$&$Y\oplus V$&$Z\oplus W$\\\hline
    $Z$&&&&$1\oplus V$&$Z\oplus W$&$Y\oplus V$\\\hline
    $V$&&&&&$1\oplus Y\oplus V$&$X\oplus Z\oplus W$\\\hline
    $W$&&&&&&$1\oplus Y\oplus V$
\end{tabular}.
\end{center}
\end{table}
\hspace{-17pt}Thus, they have
\[ \hspace{-30pt}\fp_{\mcal B}(1)=1=\fp_{\mcal B}(X),\quad\fp_{\mcal B}(Y)=\frac{\sin\frac{2\pi}7}{\sin\frac\pi7}=\fp_{\mcal B}(Z),\quad\fp_{\mcal B}(V)=\frac{\sin\frac{3\pi}7}{\sin\frac\pi7}=\fp_{\mcal B}(W), \]
and
\[ \fp(\mc B)=\frac7{2\sin^2\frac\pi7}\approx18.6. \]
There are six solutions of quantum dimensions
\begin{align*}
    (d_X,d_Y,d_Z,d_V,d_W)=&(-1,\frac{\sin\frac\pi7}{\cos\frac\pi{14}},-\frac{\sin\frac\pi7}{\cos\frac\pi{14}},-\frac{\sin\frac{2\pi}7}{\cos\frac\pi{14}},\frac{\sin\frac{2\pi}7}{\cos\frac\pi{14}}),(1,\frac{\sin\frac\pi7}{\cos\frac\pi{14}},\frac{\sin\frac\pi7}{\cos\frac\pi{14}},-\frac{\sin\frac{2\pi}7}{\cos\frac\pi{14}},-\frac{\sin\frac{2\pi}7}{\cos\frac\pi{14}}),\\
    &(-1,-\frac{\sin\frac{3\pi}7}{\cos\frac{3\pi}{14}},\frac{\sin\frac{3\pi}7}{\cos\frac{3\pi}{14}},\frac{\sin\frac\pi7}{\cos\frac{3\pi}{14}},-\frac{\sin\frac\pi7}{\cos\frac{3\pi}{14}}),(1,-\frac{\sin\frac{3\pi}7}{\cos\frac{3\pi}{14}},-\frac{\sin\frac{3\pi}7}{\cos\frac{3\pi}{14}},\frac{\sin\frac\pi7}{\cos\frac{3\pi}{14}},\frac{\sin\frac\pi7}{\cos\frac{3\pi}{14}}),\\
    &(-1,\frac{\sin\frac{2\pi}7}{\sin\frac\pi7},-\frac{\sin\frac{2\pi}7}{\sin\frac\pi7},\frac{\sin\frac{3\pi}7}{\sin\frac\pi7},-\frac{\sin\frac{3\pi}7}{\sin\frac\pi7}),(1,\frac{\sin\frac{2\pi}7}{\sin\frac\pi7},\frac{\sin\frac{2\pi}7}{\sin\frac\pi7},\frac{\sin\frac{3\pi}7}{\sin\frac\pi7},\frac{\sin\frac{3\pi}7}{\sin\frac\pi7}).
\end{align*}
The last solution gives unitary MFCs. They have categorical dimensions
\[ D^2(\mc B)=\frac7{2\cos^2\frac\pi{14}}(\approx3.7),\quad\frac7{2\cos^2\frac{3\pi}{14}}(\approx5.7),\quad\frac7{2\sin^2\frac\pi7}. \]
They have conformal dimensions
\[ \hspace{-50pt}(h_X,h_Y,h_Z,h_V,h_W)=\begin{cases}(\frac14,\frac37,\frac{19}{28},\frac17,\frac{11}{28}),(\frac14,\frac47,\frac{23}{28},\frac67,\frac3{28}),(\frac34,\frac37,\frac5{28},\frac17,\frac{25}{28}),(\frac34,\frac47,\frac9{28},\frac67,\frac{17}{28})&(\text{1st\&2nd }d_j\text{'s}),\\(\frac14,\frac27,\frac{15}{28},\frac37,\frac{19}{28}),(\frac14,\frac57,\frac{27}{28},\frac47,\frac{23}{28}),(\frac34,\frac27,\frac1{28},\frac37,\frac5{28}),(\frac34,\frac57,\frac{13}{28},\frac47,\frac9{28})&(\text{3rd\&4th }d_j\text{'s}),\\(\frac14,\frac17,\frac{11}{28},\frac57,\frac{27}{28}),(\frac14,\frac67,\frac3{28},\frac27,\frac{15}{28}),(\frac34,\frac17,\frac{25}{28},\frac57,\frac{13}{28}),(\frac34,\frac67,\frac{17}{28},\frac27,\frac1{28})&(\text{5th\&6th }d_j\text{'s}).\end{cases}\quad(\mods1) \]
Therefore, there are
\[ 6(\text{quantum dimensions})\times4(\text{conformal dimensions})\times2(\text{categorical dimensions})=48 \]
MFC, among which those eight with the last quantum dimensions give unitary MFCs. We classify connected étale algebras in all 48 MFCs simultaneously.

There are 12 sets of natural numbers which obey (\ref{FPdimA2bound})
\begin{align*}
    (n_X,n_Y,n_Z,n_V,n_W)=&(0,0,0,0,0),(1,0,0,0,0),(2,0,0,0,0),(3,0,0,0,0),\\
    &(1,1,0,0,0),(1,0,1,0,0),(1,0,0,1,0),(1,0,0,0,1),\\
    &(0,1,0,0,0),(0,0,1,0,0),(0,0,0,1,0),(0,0,0,0,1).
\end{align*}
The nontrivial 11 solutions do not give commutative algebra because they contain object(s) with nontrivial conformal dimensions and fail to satisfy the necessary condition.

To summarize, we found
\begin{table}[H]
\begin{center}
\begin{tabular}{c|c|c|c}
    Connected étale algebra $A$&$\mc B_A$&$\rank(\mc B_A)$&Lagrangian?\\\hline
    $1$&$\mc B$&$6$&No
\end{tabular}.
\end{center}
\caption{Connected étale algebras in rank six MFC $\mcal B\simeq su(2)_5\simeq\vecG_{\mbb Z/2\mbb Z}^{-1}\boxtimes psu(2)_5$}\label{rank6su25results}
\end{table}

\subsubsection{$\mc B\simeq so(5)_2$}\label{so52}
The six simple objects obey fusion rules
\begin{table}[H]
\begin{center}
\begin{tabular}{c|c|c|c|c|c|c}
    $\otimes$&$1$&$X$&$Y$&$Z$&$V$&$W$\\\hline
    $1$&$1$&$X$&$Y$&$Z$&$V$&$W$\\\hline
    $X$&&$1$&$Y$&$Z$&$W$&$V$\\\hline
    $Y$&&&$1\oplus X\oplus Z$&$Y\oplus Z$&$V\oplus W$&$V\oplus W$\\\hline
    $Z$&&&&$1\oplus X\oplus Y$&$V\oplus W$&$V\oplus W$\\\hline
    $V$&&&&&$1\oplus Y\oplus Z$&$X\oplus Y\oplus Z$\\\hline
    $W$&&&&&&$1\oplus Y\oplus Z$
\end{tabular}.
\end{center}
\end{table}
\hspace{-17pt}Thus, they have
\[ \fp_{\mc B}(1)=1=\fp_{\mc B}(X),\quad\fp_{\mc B}(Y)=2=\fp_{\mc B}(Z),\quad\fp_{\mc B}(V)=\sqrt5=\fp_{\mc B}(W), \]
and
\[ \fp(\mc B)=20. \]
For quantum dimension, there are two solutions
\[ (d_X,d_Y,d_Z,d_V,d_W)=(1,2,2,-\sqrt5,-\sqrt5),(1,2,2,\sqrt5,\sqrt5). \]
The second solution gives unitary MFCs. They have the same categorical dimension
\[ D^2(\mc B)=20. \]
Regardless of quantum dimensions, there are four conformal dimensions\footnote{Naively, one finds 16 conformal dimensions, but the others are related to one in the main text under permutations $(YZ)$ or $(VW)$ of simple objects.}
\[ (h_X,h_Y,h_Z,h_V,h_W)=(0,\frac15,\frac45,0,\frac12),(0,\frac15,\frac45,\frac14,\frac34),(0,\frac25,\frac35,0,\frac12),(0,\frac25,\frac35,\frac14,\frac34).\quad(\mods1) \]
Therefore, there are
\[ 2(\text{quantum dimensions})\times4(\text{conformal dimensions})\times2(\text{categorical dimensions})=16 \]
MFCs, among which those eight with the second quantum dimensions give unitary MFCs. We classify connected étale algebras in all 16 MFCs simultaneously.

For this to obey (\ref{FPdimA2bound}), 12 sets of natural numbers are available
\begin{align*}
    (n_X,n_Y,n_Z,n_V,n_W)=&(0,0,0,0,0),(1,0,0,0,0),(2,0,0,0,0),(3,0,0,0,0),\\
    &(1,1,0,0,0),(1,0,1,0,0),(1,0,0,1,0),(1,0,0,0,1),\\
    &(0,1,0,0,0),(0,0,1,0,0),(0,0,0,1,0),(0,0,0,0,1).
\end{align*}
The first solution is nothing but the trivial connected étale algebra $A\cong1$ giving $\mc B_A^0\simeq\mc B_A\simeq\mc B$.

The second solution gives the $\mbb Z/2\mbb Z$ algebra $A\cong1\oplus X$. Since the $\mbb Z/2\mbb Z$ object $X$ has $(d_X,h_X)=(1,0)$, it does give connected étale algebra \cite{KK23preMFC}. Let us determine the category of right $A$-modules. The algebra has $\fp_{\mc B}(A)=2$, and demands
\[ \fp(\mc B_A^0)=5,\quad\fp(\mc B_A)=10. \]
For the first category, the only possibility is
\[ \mc B_A^0\simeq\vecG_{\mbb Z/5\mbb Z}^1. \]
(They both have additive central charges $c=0$ mod 4.) The second category contains this $\mbb Z/5\mbb Z$ MFC as a subcategory. It turns out
\[ \mc B_A\simeq\text{TY}(\mbb Z/5\mbb Z), \]
a $\mbb Z/5\mbb Z$ Tambara-Yamagami category \cite{TY98}. One of the easiest ways to find this fact is to perform anyon condensation. In the process, we `identify' $X$ with $1$. Since $W\cong X\otimes V,V\cong X\otimes W$, we further `identify' $V$ and $W$. Invariant $Y,Z$ with quantum dimensions two split into two each. As a result, we get five invertible simple objects, and one simple object with quantum dimension $\pm\sqrt5$. This is nothing but the $\mbb Z/5\mbb Z$ Tambara-Yamagami category.\footnote{More rigorously, we have to find NIM-reps. Indeed, we find a six-dimensional NIM-rep
\[ \hspace{-20pt}n_1=1_6=n_X,\quad n_Y=\begin{pmatrix}0&0&0&1&1&0\\0&0&1&0&1&0\\0&1&0&1&0&0\\1&0&1&0&0&0\\1&1&0&0&0&0\\0&0&0&0&0&2\end{pmatrix},\quad n_Z=\begin{pmatrix}0&1&1&0&0&0\\1&0&0&1&0&0\\1&0&0&0&1&0\\0&1&0&0&1&0\\0&0&1&1&0&0\\0&0&0&0&0&2\end{pmatrix},\quad n_V=\begin{pmatrix}0&0&0&0&0&1\\0&0&0&0&0&1\\0&0&0&0&0&1\\0&0&0&0&0&1\\0&0&0&0&0&1\\1&1&1&1&1&0\end{pmatrix}=n_W .\]
Denoting a basis of $\mc B_A$ by $\{m_1,m_2,m_3,m_4,m_5,m_6\}$, we get a multiplication table
\begin{table}[H]
\begin{center}
\makebox[1 \textwidth][c]{       
\resizebox{1.2 \textwidth}{!}{\begin{tabular}{c|c|c|c|c|c|c}
    $b_j\otimes\backslash$&$F(m_1)$&$F(m_2)$&$F(m_3)$&$F(m_4)$&$F(m_5)$&$F(m_6)$\\\hline
    1&$F(m_1)$&$F(m_2)$&$F(m_3)$&$F(m_4)$&$F(m_5)$&$F(m_6)$\\
    $X$&$F(m_1)$&$F(m_2)$&$F(m_3)$&$F(m_4)$&$F(m_5)$&$F(m_6)$\\
    $Y$&$F(m_4)\oplus F(m_5)$&$F(m_3)\oplus F(m_5)$&$F(m_2)\oplus F(m_4)$&$F(m_1)\oplus F(m_3)$&$F(m_1)\oplus F(m_2)$&$2F(m_6)$\\
    $Z$&$F(m_2)\oplus F(m_3)$&$F(m_1)\oplus F(m_4)$&$F(m_1)\oplus F(m_5)$&$F(m_2)\oplus F(m_5)$&$F(m_3)\oplus F(m_4)$&$2F(m_6)$\\
    $V$&$F(m_6)$&$F(m_6)$&$F(m_6)$&$F(m_6)$&$F(m_6)$&$F(m_1)\oplus F(m_2)\oplus F(m_3)\oplus F(m_4)\oplus F(m_5)$\\
    $W$&$F(m_6)$&$F(m_6)$&$F(m_6)$&$F(m_6)$&$F(m_6)$&$F(m_1)\oplus F(m_2)\oplus F(m_3)\oplus F(m_4)\oplus F(m_5)$
\end{tabular}.}}
\end{center}
\end{table}
\hspace{-14pt}In the basis, the monoidal products give the identifications
\[ F(m_1)\cong1\oplus X,\quad F(m_2)\cong Z\cong F(m_3),\quad F(m_4)\cong Y\cong F(m_5),\quad F(m_6)\cong V\oplus W. \]
In the category $\mc B_A$ of right $A$-modules, they have quantum dimensions \cite{KO01}
\[ d_{\mc B_A}(F(m_1))=d_{\mc B_A}(F(m_2))=d_{\mc B_A}(F(m_3))=d_{\mc B_A}(F(m_4))=d_{\mc B_A}(F(m_5))=1,\quad d_{\mc B_A}(F(m_6))=\pm\sqrt5, \]
showing $\mc B_A\simeq\text{TY}(\mbb Z/5\mbb Z)$, especially $\rank(\mc B_A)=6$.}

The third and fourth solutions have $\fp_{\mc B}(A)=3,4$, and demands $\fp(\mc B_A^0)=\frac{20}9,\frac54$, but there is no MFC with these Frobenius-Perron dimension. Thus, the two candidates are ruled out. Other four solutions $A\cong1\oplus X\oplus V,1\oplus X\oplus W,1\oplus V,1\oplus W$ have $\fp_{\mc B}=2+\sqrt5,1+\sqrt5$, and they are ruled out for the same reason.

The other four solutions with $Y$ or $Z$ fail to be commutative because they have nontrivial conformal dimensions. Thus, they are also ruled out.

To sum up, we found
\begin{table}[H]
\begin{center}
\begin{tabular}{c|c|c|c}
    Connected étale algebra $A$&$\mc B_A$&$\rank(\mc B_A)$&Lagrangian?\\\hline
    $1$&$\mc B$&$6$&No\\
    $1\oplus X$&$\text{TY}(\mbb Z/5\mbb Z)$&6&No
\end{tabular}.
\end{center}
\caption{Connected étale algebras in rank six MFC $\mcal B\simeq so(5)_2$}\label{rank6so52results}
\end{table}
\hspace{-17pt}Namely, all 16 MFCs $\mc B\simeq so(5)_2$'s fail to be completely anisotropic.

\subsubsection{$\mc B\simeq\fib\boxtimes psu(2)_5$}
We have the following fusion table
\begin{table}[H]
\begin{center}
\begin{tabular}{c|c|c|c|c|c|c}
    $\otimes$&$1$&$X$&$Y$&$Z$&$V$&$W$\\\hline
    $1$&$1$&$X$&$Y$&$Z$&$V$&$W$\\\hline
    $X$&&$1\oplus X$&$V$&$W$&$Y\oplus V$&$Z\oplus W$\\\hline
    $Y$&&&$1\oplus Z$&$Y\oplus Z$&$X\oplus W$&$V\oplus W$\\\hline
    $Z$&&&&$1\oplus Y\oplus Z$&$V\oplus W$&$X\oplus V\oplus W$\\\hline
    $V$&&&&&$1\oplus X\oplus Z\oplus W$&$Y\oplus Z\oplus V\oplus W$\\\hline
    $W$&&&&&&$1\oplus X\oplus Y\oplus Z\oplus V\oplus W$
\end{tabular}.
\end{center}
\end{table}
\hspace{-17pt}The Frobenius-Perron dimension is given by
\begin{align*}
    \fp_{\mc B}(1)=1,\quad\fp_{\mc B}(X)=&\zeta,\quad\fp_{\mc B}(Y)=\frac{\sin\frac{5\pi}7}{\sin\frac\pi7},\\
    \fp_{\mc B}(Z)=\frac{\sin\frac{3\pi}7}{\sin\frac\pi7},\quad\fp_{\mc B}(V)=&\zeta\frac{\sin\frac{5\pi}7}{\sin\frac\pi7},\quad\fp_{\mc B}(W)=\zeta\frac{\sin\frac{3\pi}7}{\sin\frac\pi7},
\end{align*}
and
\[ \fp(\mc B)=\frac{5+\sqrt5}2\frac7{4\sin^2\frac\pi7}\approx33.6. \]
Since the MFCs are given by Deligne tensor products of two factors, the two factors should be modular in order to get modular $\mc B$. This is automatic for the two factors $\fib,psu(2)_5$. Thanks to the product structure, we  know their quantum
\begin{align*}
    \hspace{-57pt}(d_X,d_Y,d_Z,d_V,d_W)=&(-\zeta^{-1},\frac{\sin\frac\pi7}{\cos\frac\pi{14}},-\frac{\sin\frac{2\pi}7}{\cos\frac\pi{14}},-\zeta^{-1}\frac{\sin\frac\pi7}{\cos\frac\pi{14}},\zeta^{-1}\frac{\sin\frac{2\pi}7}{\cos\frac\pi{14}}),(-\zeta^{-1},-\frac{\sin\frac{3\pi}7}{\cos\frac{3\pi}{14}},\frac{\sin\frac\pi7}{\cos\frac{3\pi}{14}},\zeta^{-1}\frac{\sin\frac{3\pi}7}{\cos\frac{3\pi}{14}},-\zeta^{-1}\frac{\sin\frac\pi7}{\cos\frac{3\pi}{14}}),\\
    &(-\zeta^{-1},\frac{\sin\frac{5\pi}7}{\sin\frac\pi7},\frac{\sin\frac{3\pi}7}{\sin\frac\pi7},-\zeta^{-1}\frac{\sin\frac{5\pi}7}{\sin\frac\pi7},-\zeta^{-1}\frac{\sin\frac{3\pi}7}{\sin\frac\pi7}),(\zeta,\frac{\sin\frac\pi7}{\cos\frac\pi{14}},-\frac{\sin\frac{2\pi}7}{\cos\frac\pi{14}},\zeta\frac{\sin\frac\pi7}{\cos\frac\pi{14}},-\zeta\frac{\sin\frac{2\pi}7}{\cos\frac\pi{14}}),\\
    &(\zeta,-\frac{\sin\frac{3\pi}7}{\cos\frac{3\pi}{14}},\frac{\sin\frac\pi7}{\cos\frac{3\pi}{14}},-\zeta\frac{\sin\frac{3\pi}7}{\cos\frac{3\pi}{14}},\zeta\frac{\sin\frac\pi7}{\cos\frac{3\pi}{14}}),(\zeta,\frac{\sin\frac{5\pi}7}{\sin\frac\pi7},\frac{\sin\frac{3\pi}7}{\sin\frac\pi7},\zeta\frac{\sin\frac{5\pi}7}{\sin\frac\pi7},\zeta\frac{\sin\frac{3\pi}7}{\sin\frac\pi7}),
\end{align*}
categorical
\begin{align*}
    D^2(\mc B)=&\frac{5-\sqrt5}2\frac7{4\cos^2\frac\pi{14}}(\approx2.5),\quad\frac{5-\sqrt5}2\frac7{4\cos^2\frac{3\pi}{14}}(\approx4.0),\quad\frac{5-\sqrt5}2\frac7{4\sin^2\frac\pi7}(\approx12.8),\\
    &\frac{5+\sqrt5}2\frac7{4\cos^2\frac\pi{14}}(\approx6.7),\quad\frac{5+\sqrt5}2\frac7{4\cos^2\frac{3\pi}{14}}(\approx10.4),\quad\frac{5+\sqrt5}2\frac7{4\sin^2\frac\pi7},
\end{align*}
and conformal dimensions
\[ \hspace{-40pt}(h_X,h_Y,h_Z,h_V,h_W)=\begin{cases}(\frac15,\frac37,\frac17,\frac{22}{35},\frac{12}{35}),(\frac15,\frac47,\frac67,\frac{27}{35},\frac2{35}),(\frac45,\frac37,\frac17,\frac8{35},\frac{33}{35}),(\frac45,\frac47,\frac67,\frac{13}{35},\frac{23}{35})&(\text{1st }d_j\text{'s}),\\(\frac15,\frac27,\frac37,\frac{17}{35},\frac{22}{35}),(\frac15,\frac57,\frac47,\frac{32}{35},\frac{27}{35}),(\frac45,\frac27,\frac37,\frac3{35},\frac8{35}),(\frac45,\frac57,\frac47,\frac{18}{35},\frac{13}{35})&(\text{2nd }d_j\text{'s}),\\(\frac15,\frac17,\frac57,\frac{12}{35},\frac{32}{35}),(\frac15,\frac67,\frac27,\frac2{35},\frac{17}{35}),(\frac45,\frac17,\frac57,\frac{33}{35},\frac{18}{35}),(\frac45,\frac67,\frac27,\frac{23}{35},\frac3{35})&(\text{3rd }d_j\text{'s}),\\(\frac25,\frac37,\frac17,\frac{29}{35},\frac{19}{35}),(\frac25,\frac47,\frac67,\frac{34}{35},\frac9{35}),(\frac35,\frac37,\frac17,\frac1{35},\frac{26}{35}),(\frac35,\frac47,\frac67,\frac6{35},\frac{16}{35})&(\text{4th }d_j\text{'s}),\\(\frac25,\frac27,\frac37,\frac{24}{35},\frac{29}{35}),(\frac25,\frac57,\frac47,\frac4{35},\frac{34}{35}),(\frac35,\frac27,\frac37,\frac{31}{35},\frac1{35}),(\frac35,\frac57,\frac47,\frac{11}{35},\frac6{35})&(\text{5th }d_j\text{'s}),\\(\frac25,\frac17,\frac57,\frac{19}{35},\frac4{35}),(\frac25,\frac67,\frac27,\frac9{35},\frac{24}{35}),(\frac35,\frac17,\frac57,\frac{26}{35},\frac{11}{35}),(\frac35,\frac67,\frac27,\frac{16}{35},\frac{31}{35})&(\text{6th }d_j\text{'s}).\end{cases}\quad(\mods1) \]
Therefore, there are
\[ 6(\text{quantum dimensions})\times4(\text{conformal dimensions})\times2(\text{categorical dimensions})=48 \]
MFC, among which those eight with the last quatum dimensions give unitary MFCs. We classify connected étale algebras in all 48 MFCs simultaneously.

For this to obey (\ref{FPdimA2bound}), the natural numbers can take 94 sets. However, the nontrivial conformal dimensions of each simple objects left us with only 
the trivial connected étale algebra giving $\mc B_A^0\simeq\mc B_A\simeq\mc B$.

To summarize, we found
\begin{table}[H]
\begin{center}
\begin{tabular}{c|c|c|c}
    Connected étale algebra $A$&$\mc B_A$&$\rank(\mc B_A)$&Lagrangian?\\\hline
    $1$&$\mc B$&$6$&No
\end{tabular}.
\end{center}
\caption{Connected étale algebras in rank six MFC $\mcal B\simeq\fib\boxtimes psu(2)_5$}\label{rank6fibpsu25results}
\end{table}

\subsubsection{$\mc B\simeq psu(2)_{11}$}
The fusion rules of this MFC are
\begin{table}[H]
\begin{center}
\begin{tabular}{c|c|c|c|c|c|c}
    $\otimes$&$1$&$X$&$Y$&$Z$&$V$&$W$\\\hline
    $1$&$1$&$X$&$Y$&$Z$&$V$&$W$\\\hline
    $X$&&$1\oplus Y$&$X\oplus Z$&$Y\oplus V$&$Z\oplus W$&$V\oplus W$\\\hline
    $Y$&&&$1\oplus Y\oplus V$&$X\oplus Z\oplus W$&$Y\oplus V\oplus W$&$Z\oplus V\oplus W$\\\hline
    $Z$&&&&$1\oplus Y\oplus V\oplus W$&$X\oplus Z\oplus V\oplus W$&$Y\oplus Z\oplus V\oplus W$\\\hline
    $V$&&&&&$1\oplus Y\oplus Z\oplus V\oplus W$&$X\oplus Y\oplus Z\oplus V\oplus W$\\\hline
    $W$&&&&&&$1\oplus X\oplus Y\oplus Z\oplus V\oplus W$
\end{tabular}.
\end{center}
\end{table}
\hspace{-17pt}Thus, they have
\begin{align*}
    \fp_{\mc B}(1)=1,\quad\fp_{\mc B}(X)&=\frac{\sin\frac{11\pi}{13}}{\sin\frac\pi{13}},\quad\fp_{\mc B}(Y)=\frac{\sin\frac{3\pi}{13}}{\sin\frac\pi{13}},\\
    \fp_{\mc B}(Z)=\frac{\sin\frac{9\pi}{13}}{\sin\frac\pi{13}},\quad\fp_{\mc B}(V)&=\frac{\sin\frac{5\pi}{13}}{\sin\frac\pi{13}},\quad\fp_{\mc B}(W)=\frac{\sin\frac{7\pi}{13}}{\sin\frac\pi{13}},
\end{align*}
and
\[ \fp(\mc B)=\frac{13}{4\sin^2\frac\pi{13}}\approx56.7. \]
There are six solutions of quantum dimensions
\begin{align*}
    \hspace{-45pt}(d_X,d_Y,d_Z,d_V,d_W)=&(-\frac{\sin\frac\pi{13}}{\cos\frac\pi{26}},-\frac{\sin\frac{5\pi}{13}}{\cos\frac\pi{26}},\frac{\sin\frac{2\pi}{13}}{\cos\frac\pi{26}},\frac{\sin\frac{4\pi}{13}}{\cos\frac\pi{26}},-\frac{\sin\frac{3\pi}{13}}{\cos\frac\pi{26}}),(\frac{\sin\frac{3\pi}{13}}{\cos\frac{3\pi}{26}},-\frac{\sin\frac{2\pi}{13}}{\cos\frac{3\pi}{26}},-\frac{\sin\frac{6\pi}{13}}{\cos\frac{3\pi}{26}},-\frac{\sin\frac{\pi}{13}}{\cos\frac{3\pi}{26}},\frac{\sin\frac{4\pi}{13}}{\cos\frac{3\pi}{26}}),\\&(-\frac{\sin\frac{5\pi}{13}}{\cos\frac{5\pi}{26}},\frac{\sin\frac{\pi}{13}}{\cos\frac{5\pi}{26}},\frac{\sin\frac{3\pi}{13}}{\cos\frac{5\pi}{26}},-\frac{\sin\frac{6\pi}{13}}{\cos\frac{5\pi}{26}},\frac{\sin\frac{2\pi}{13}}{\cos\frac{5\pi}{26}}),(\frac{\sin\frac{6\pi}{13}}{\cos\frac{7\pi}{26}},\frac{\sin\frac{4\pi}{13}}{\cos\frac{7\pi}{26}},\frac{\sin\frac{\pi}{13}}{\cos\frac{7\pi}{26}},-\frac{\sin\frac{2\pi}{13}}{\cos\frac{7\pi}{26}},-\frac{\sin\frac{5\pi}{13}}{\cos\frac{7\pi}{26}}),\\
    &(-\frac{\sin\frac{4\pi}{13}}{\cos\frac{9\pi}{26}},\frac{\sin\frac{6\pi}{13}}{\cos\frac{9\pi}{26}},-\frac{\sin\frac{5\pi}{13}}{\cos\frac{9\pi}{26}},\frac{\sin\frac{3\pi}{13}}{\cos\frac{9\pi}{26}},-\frac{\sin\frac{\pi}{13}}{\cos\frac{9\pi}{26}}),(\frac{\sin\frac{11\pi}{13}}{\sin\frac\pi{13}},\frac{\sin\frac{3\pi}{13}}{\sin\frac\pi{13}},\frac{\sin\frac{9\pi}{13}}{\sin\frac\pi{13}},\frac{\sin\frac{5\pi}{13}}{\sin\frac\pi{13}},\frac{\sin\frac{7\pi}{13}}{\sin\frac\pi{13}}).
\end{align*}
They have categorical dimensions
\[ \hspace{-55pt}D^2(\mc B)=\frac{13}{4\cos^2\frac\pi{26}}(\approx3.3),\quad\frac{13}{4\cos^2\frac{3\pi}{26}}(\approx3.7),\quad\frac{13}{4\cos^2\frac{5\pi}{26}}(\approx4.8),\quad\frac{13}{4\cos^2\frac{7\pi}{26}}(\approx7.4),\quad\frac{13}{4\cos^2\frac{9\pi}{26}}(\approx15.0),\quad\frac{13}{4\sin^2\frac\pi{13}}, \]
respectively. They have conformal dimensions
\[ (h_X,h_Y,h_Z,h_V,h_W)=\begin{cases}(\frac2{13},\frac1{13},\frac{10}{13},\frac3{13},\frac6{13}),(\frac{11}{13},\frac{12}{13},\frac3{13},\frac{10}{13},\frac7{13})&(\text{1st quantum dimensions}),\\(\frac6{13},\frac3{13},\frac4{13},\frac9{13},\frac5{13}),(\frac7{13},\frac{10}{13},\frac9{13},\frac4{13},\frac8{13})&(\text{2nd quantum dimensions}),\\(\frac3{13},\frac8{13},\frac2{13},\frac{11}{13},\frac9{13}),(\frac{10}{13},\frac5{13},\frac{11}{13},\frac2{13},\frac4{13})&(\text{3rd quantum dimensions}),\\(\frac1{13},\frac7{13},\frac5{13},\frac8{13},\frac3{13}),(\frac{12}{13},\frac6{13},\frac8{13},\frac5{13},\frac{10}{13})&(\text{4th quantum dimensions}),\\(\frac5{13},\frac9{13},\frac{12}{13},\frac1{13},\frac2{13}),(\frac8{13},\frac4{13},\frac1{13},\frac{12}{13},\frac{11}{13})&(\text{5th quantum dimensions}),\\(\frac4{13},\frac2{13},\frac7{13},\frac6{13},\frac{12}{13}),(\frac9{13},\frac{11}{13},\frac6{13},\frac7{13},\frac1{13})&(\text{6th quantum dimensions}).\end{cases}\quad(\mods1) \]
Therefore, there are
\[ 6(\text{quantum dimensions})\times2(\text{conformal dimensions})\times2(\text{categorical dimensions})=24 \]
MFCs, among which those four with the last quantum dimensions give unitary MFCs. We classify connected étale algebras in all 24 MFCs simultaneously.

Then, the natural numbers $n_j$'s can take only 14 values
\begin{align*}
    (n_X,n_Y,n_Z,n_V,n_W)=&(0,0,0,0,0),(1,0,0,0,0),(2,0,0,0,0),\\
    &(3,0,0,0,0),(1,1,0,0,0),(1,0,1,0,0),\\
    &(1,0,0,1,0),(1,0,0,0,1),(0,1,0,0,0),\\
    &(0,2,0,0,0),(0,1,1,0,0),(0,0,1,0,0),\\
    &(0,0,0,1,1),(0,0,0,0,1).
\end{align*}
The nontrivial 13 candidates contain simple object(s) $b_j$ with nontrivial conformal dimensions, and do not give connected étale algebra.

We conclude
\begin{table}[H]
\begin{center}
\begin{tabular}{c|c|c|c}
    Connected étale algebra $A$&$\mc B_A$&$\rank(\mc B_A)$&Lagrangian?\\\hline
    $1$&$\mc B$&$6$&No
\end{tabular}.
\end{center}
\caption{Connected étale algebras in rank six MFC $\mcal B\simeq psu(2)_{11}$}\label{rank6fibpsu211results}
\end{table}

\section{Physical applications}
\subsection{Theorems}
In this section, we discuss the physical applications of the classification results.

Consider a two-dimensional $\mc C$-symmetric gapped phase. Ground state degeneracy (GSD) is related to module category. It is known \cite{TW19,HLS21} that 
\[
\{\text{2d }\mathcal{C}\text{-symmetric gapped phases}\}\cong\{\mathcal{C}\text{{-module categories} }\mathcal{M}\}.
\]
After the $\mc C$ symmetry is spontaneously broken, the physical system will create several ground states. The ground state degeneracy (in the LHS) is determined by the rank of module categories $\mathcal{M}$ (in the RHS). Therefore, the physical problems in the LHS can be translated to mathematical problems in the RHS. In particular, ground state degeneracy is given by
\begin{center}
    GSD\ =\ rank($\mathcal{M}$).
\end{center}
This leads to the\\

\textbf{Theorem}.\ \emph{Let $\mc B$ be a rank six multiplicity-free modular fusion category and $A \in\mc B$ be a connected étale algebra. Suppose two-dimensional $\mc B$-symmetric grapped phases are described by indecomposable $\mc B_A$'s. Then, the gapped phases have}
\begin{align*}
\text{GSD}\in
\begin{cases}
\{6\}&(\mc B\simeq \vecG
_{\mbb Z/6\mbb Z}^\alpha),\\
\{6\}&(\mc B\simeq\vecG_{\mbb Z/2\mbb Z}^{-1}\boxtimes\ising,\\
\{6\}&(\mc B\simeq su(3)_2\simeq\fib\boxtimes\vecG_{\mbb Z/3\mbb Z}^1),\\
\{6\}&(\mc B\simeq\text{TriCritIsing}),\\
\{6\}&(\mc B\simeq su(2)_5\simeq\vecG_{\mbb Z/2\mbb Z}^{-1}\boxtimes psu(2)_5),\\
\{6\}&(\mc B\simeq so(5)_2),\\
\{6\}&(\mc B\simeq\fib\boxtimes psu(2)_{5}),\\
\{6\}&(\mc B\simeq psu(2)_{11}).
\end{cases}
\end{align*}\newline

The Theorem also proves certain SSBs. Here, we have the\\

\textbf{Definition}. \cite{KK23GSD} Let $\mc C$ be a fusion category and  $\mc M$ a (left) $\mc C$-module category describing a $\mc C$-symmetric gapped phase. For a symmetry $c\in\mc C$, if $\exists m\in \mc M$ such that $c\triangleright m\ncong m$, then we call \emph{c} is \emph{spontaneously broken}. We also say $\mc C$ is \emph{spontaneously broken} if there is a spontaneously broken object $\emph{c}\in\mc C$. Otherwise, the category symmetry $\mc C$ is called \emph{preserved}, namely all objects act trivially.\\

With the definition, one can show a\\

\textbf{Lemma}. \cite{KK23GSD} \emph{Let $\mc C$ be a fusion category and $\mc M$ be an indecomposable (left) $\mc C$-module category. Then $rank(\mc M)>1$ implies $\mc C$ is spontaneously broken.}\\

Therefore, we have proved SSBs:\\

\textbf{Theorem}. \emph{Let $\mc B$ be a modular fusion category, and $A \in\mc B $ be a connected étale algebra. In the two-dimensional gapped phases decribled by $\mc B_A$'s, $\mc B$ symmetries are spontaneously broken for}
\begin{align*}
\mc B\simeq
\begin{cases}
\vecG_{\mbb Z/6\mbb Z}^\alpha,\\
\vecG_{\mbb Z/2\mbb Z}^{-1}\boxtimes\ising,\\
su(3)_2\simeq\fib\boxtimes\vecG_{\mbb Z/3\mbb Z}^1,\\
\text{TriCritIsing},\\
\vecG_{\mbb Z/2\mbb Z}^{-1}\boxtimes psu(2)_5,\\
so(5)_2,\\
\fib\boxtimes psu(2)_{5},\\
psu(2)_{11}.
\end{cases}
\end{align*}\newline

\textbf{Remark}. As noted in \cite{KK23GSD}, commutativity of connected étale algebra seems too strong; numerical computation suggests an existence of $\mc B$-symmetric gapped phase described by $\mc B_A$ with non-commutative
connected separable algebra $A\in\mc B$.

\subsection{Examples}
In this section, we discuss concrete examples and \textit{predict} GSDs and SSB. As in \cite{KK23rank5}, we consider relevant deformations of non-unitary minimal models.

Pick a non-unitary minimal model\footnote{We basically follow the notations of \cite{FMS}.} $M(p,2p\pm1)$ with $p\ge2$. It was proved \cite{KK22free} that the relevant $\phi_{5,1}$-deformation of $M(p,2p+1)$ preserves rank $(p-1)$ MFC formed by Verlinde lines $\{\mc L_{1,1},\mc L_{1,2},\dots,\mc L_{1,p-1}\}$, and the relevant $\phi_{1,2}$-deformation of $M(p,2p-1)$ preserves rank $(p-1)$ MFC formed by $\{\mc L_{1,1},\mc L_{3,1},\dots,\mc L_{2p-3,1}\}$. For $p=7$, the preserved MFCs have rank six, and would fall in our classifications (if they are multiplicity-free). Let us study the two examples in more detail below.\newline

\paragraph{$M(7,15)+\phi_{5,1}$.} The relevant deformation preserves rank six MFC with simple objects $\{\mc L_{1,1},\mc L_{1,2},\mc L_{1,3},\mc L_{1,4},\mc L_{1,5},\mc L_{1,6}\}$. They form $\mc B\simeq su(2)_5$ with identifications
\[ 1\cong\mc L_{1,1},\quad X\cong\mc L_{1,6},\quad Y\cong\mc L_{1,5},\quad Z\cong\mc L_{1,2},\quad V\cong\mc L_{1,3},\quad W\cong\mc L_{1,4}. \]
One can read this off from their (non-unitary) quantum dimensions
\[ (d_{1,1},d_{1,6},d_{1,5},d_{1,2},d_{1,3},d_{1,4})=(1,-1,\frac{\sin\frac{2\pi}7}{\sin\frac\pi7},-\frac{\sin\frac{2\pi}7}{\sin\frac\pi7},\frac{\sin\frac{3\pi}7}{\sin\frac\pi7},-\frac{\sin\frac{3\pi}7}{\sin\frac\pi7}) \]
and conformal dimensions
\[ (h_{1,1},h_{1,6},h_{1,5},h_{1,2},h_{1,3},h_{1,4})=(0,\frac{65}4,\frac{76}7,\frac{31}{28},\frac{23}7,\frac{183}{28}). \]
One can check these coincide with our fifth quantum dimensions and its second conformal dimensions (mod 1). We also checked the fusion ring is multiplicity-free and coincides with that of $su(2)_5$. Therefore, our classification result implies the massive RG flow described by $\mc B_A$ should have $\text{GSD}=6$ and $\mc B$ symmetry should be spontaneously broken.
\newline

\paragraph{$M(7,13)+\phi_{1,2}$.} The relevant deformation preserves rank six MFC with simple objects $\{\mc L_{1,1},\mc L_{3,1},\mc L_{5,1},\mc L_{7,1},\mc L_{9,1},\mc L_{11,1}\}$. They form $\mc B\simeq psu(2)_{11}$ with identifications
\[ 1\cong\mc L_{1,1},\quad X\cong\mc L_{11,1},\quad Y\cong\mc L_{3,1},\quad Z\cong\mc L_{9,1},\quad V\cong\mc L_{5,1},\quad W\cong\mc L_{7,1}. \]
This is because they have our first (non-unitary) quantum dimensions
\[ (d_{1,1},d_{11,1},d_{3,1},d_{9,1},d_{5,1},d_{7,1})=(1,-\frac{\sin\frac{\pi}{13}}{\cos\frac\pi{26}},-\frac{\sin\frac{5\pi}{13}}{\cos\frac\pi{26}},\frac{\sin\frac{2\pi}{13}}{\cos\frac\pi{26}},\frac{\sin\frac{4\pi}{13}}{\cos\frac\pi{26}},-\frac{\sin\frac{3\pi}{13}}{\cos\frac\pi{26}}) \]
and conformal dimensions
\[ (h_{1,1},h_{11,1},h_{3,1},h_{9,1},h_{5,1},h_{7,1})=(0,\frac{145}{13},\frac1{13},\frac{88}{13},\frac{16}{13},\frac{45}{13}). \]
One sees these coincide with our first conformal dimensions mod 1. We also checked the fusion ring is multiplicity-free and coincides with that of $psu(2)_{11}$. Therefore, our classification result implies the massive RG flow described by $\mc B_A$ should have $\text{GSD}=6$ and $\mc B$ symmetry should be spontaneously broken.

\appendix
\setcounter{section}{0}
\renewcommand{\thesection}{\Alph{section}}
\setcounter{equation}{0}
\renewcommand{\theequation}{\Alph{section}.\arabic{equation}}

\end{document}